\documentclass[11pt]{amsart}
\def\from{From}

\newtheorem{theorem}{Theorem}[section]
\newtheorem{lemma}{Lemma}[section]
\newtheorem{corollary}{Corollary}[section]
\theoremstyle{definition}
\newtheorem{remark}{Remark}[section]

\numberwithin{equation}{section}

\newcommand{\NN}{{\mathbb N}}

\renewcommand{\Re}{{\rm Re}}
\newcommand{\Un}{{\mathbf 1}}
\let\wh=\widehat

\begin{document}

\title[An estimate of the zeta function]{On an analytic estimate in the
theory of the Riemann Zeta function and a Theorem of B\'aez-Duarte}

\author{Jean-Fran\c cois Burnol}

\date{18 February 2002}

\email{burnol@math.unice.fr}

\keywords{Riemann hypothesis, Nyman-Beurling theorem}

\begin{abstract}
We establish a uniform upper estimate for the values of
$\zeta(s)/\zeta(s+A)$, $0\leq A$, on the critical line (conditionally on
the Riemann Hypothesis). We use this to give a variant, purely complex
analytic, to B\'aez-Duarte's proof of a strengthened Nyman-Beurling
criterion for the validity of the Riemann Hypothesis.
\end{abstract}

\maketitle
%\parindent=0pt
%\parskip=\the\baselineskip

\section{Introduction}

The following theorem has been established by B\'aez-Duarte
\cite{luisnatural} ($\{x\} = x - [x]$ is the fractional part of the real
number $x$):

\begin{theorem}[B\'aez-Duarte \cite{luisnatural}]\label{luisthmA}
If the Riemann Hypothesis holds then the function $\Un_{0<t\leq1}$ belongs
to the closure in $L^2((0,\infty),dt)$ of the finite linear combinations of
the functions $t\to \{1/nt\}$, $n\geq1$, $n\in\NN$.
\end{theorem}

That the converse holds is (a special case of) the easy half of the
classical Nyman-Beurling criterion \cite{nym50} \cite{beur55} (in a minor
variant as the original formulation is entirely inside $L^2((0,1),dt)$.)
The not-so-easy other half of the original version of this criterion states
that the Riemann Hypothesis implies that $\Un_{0<t\leq1}$ belongs to the
span of the functions $t\to \{1/\Lambda t\}$ where $\Lambda$ varies in the
continuous interval $[1,\infty)$. The connection with the Riemann zeta
function is established with the help of the classical formula
$$\int_0^\infty \left\{\frac1t\right\} t^{s-1}\,dt = \frac{-\zeta(s)}s$$
where the integral is absolutely convergent for $0<\Re(s)<1$. Generally
speaking the Mellin transform $f(t)\mapsto \wh{f}(s)= \int_0^\infty
f(t)t^{s-1}\,dt$ establishes an isometry between $L^2((0,\infty),dt)$ and
the Hilbert space of square\--integrable functions on the critical line for
the measure $|ds|/2\pi$. So the Nyman-Beurling criterion is that the
Riemann Hypothesis holds if and only if one can approximate in the
square-mean sense the function $1/s$ on the critical line by expressions
$\zeta(s)P(s)/s$ where $P(s) = \sum_k c_k \Lambda_k^{-s}$, $\Lambda_k\geq1$.

This topic has attracted some interest in recent years in a number of
publications, among them \cite{luis93}, \cite{luis99}, \cite{luisarithm},
 \cite{balnotes3}, \cite{balnotes1}, \cite{balexpo}, \cite{con00},
\cite{lee96}, \cite{nik95}, \cite{vas96}, \cite{vas99}  and it was expected from the
expression $1/\zeta(s) = \sum_{n=1}^\infty \mu(n)n^{-s}$ that the values
$\Lambda_k = n$, $n\in \NN\setminus\{0\}$ are enough. This is what has been
established by  B\'aez-Duarte. We state again here for future reference his
theorem in complex analytic form:

\begin{theorem}[B\'aez-Duarte \cite{luisnatural}]\label{luisthmB}
If the Riemann Hypothesis holds then $1/s$ can be arbitrarily well
approximated in square-mean on the critical line by functions
$\zeta(s)P(s)/s$ with $P(s) = \sum_{n=1}^N c_{n,N}\;n^{-s}$ a Dirichlet
polynomial.
\end{theorem}

\begin{remark}
It remains an open problem to exhibit for $N\to\infty$ an explicit
sequence of ``natural approximations'' $\sum_{n=1}^N
c_{n,N}\;n^{-s}$. It is known \cite{balnotes3} that the Hilbert space
distance from $\zeta(s)P(s)/s$ to $1/s$ is bounded below
asymptotically by $C/\sqrt{\log(N)}$ for a $C>0$ (and from
\cite{jfadvances} it is known that $C\geq 2+\gamma-\log(4\pi)$.)  This
applies even to  Dirichlet polynomials $\sum_k c_k \Lambda_k^{-s}$
allowing non-integer $\Lambda_k$'s ($1\leq\Lambda_k\leq N$). It is not
known whether the restriction to integer $\Lambda$'s may result in a
worse rate of convergence. It is naturally expected that this is not
the case.  But no explicit upper bound to the Hilbert space distance
to the optimal natural approximation of degree $N$ is currently known.
\end{remark}

It is a quite notable feature of B\'aez-Duarte's proof that it frees
completely the Nyman-Beurling criterion from any appeal to the theory of
invariant subspaces of Hardy spaces. In this manner the topic of the
Nyman-Beurling criterion becomes less foreign to the more classical
analytic number theoretical topics as they are treated in Titchmarsh's book
\cite{titch51}.

The purpose of this paper is to reinforce this with the help of a purely
complex analytical proof of \ref{luisthmB}. While inspired by the original
proof and quite directly related to it, our approach relies on a novel
analytic estimate on the Riemann zeta function, conditional on the Riemann
Hypothesis:

\begin{theorem}\label{jfzeta}
Let $\epsilon>0$. Conditional on the Riemann Hypothesis one has:
$$\left|\frac{\zeta(s)}{\zeta(s+A)}\right| =
O_\epsilon(|s|^{\inf(\epsilon,A/2)})$$
on the critical line, uniformly with respect to $0\leq A<\infty$.
\end{theorem}

The other main component of the proof, as in B\'aez-Duarte's, is the use of
a Theorem of Balazard and Saias \cite[Lemme 2]{balnotes1}, which is a
version with a fine error estimate of the convergence of $\sum_{n=1}^\infty
\mu(n)/n^s$ to $1/\zeta(s)$ for $\Re(s)>1/2$. We state a slightly weakened
form of this Theorem which is directly suitable for our later use:

\begin{theorem}[Balazard-Saias \cite{balnotes1}]\label{bslemme}
Let $\frac12>\epsilon>0$ and $\theta>0$. Conditional on the Riemann
Hypothesis one has:
$$\sum_{n=1}^N \frac{\mu(n)}{n^s} = \frac1{\zeta(s)} + O_{\epsilon,
\theta}(\frac{|s|^\theta}{N^{\epsilon/3}})$$
on $\Re(s) = 1/2 + \epsilon$.
\end{theorem}

It will be apparent that the proof of the strengthened Nyman-Beurling
criterion in fact would go through equally well with much coarser error
estimates than the ones provided by \ref{jfzeta} and \ref{bslemme}. But the
future developments could possibly use the finer estimates.

This paper is organized as follows: in the first section we prove the
B\'aez-Duarte strengthened ``only  if'' Nyman-Beurling criterion using
\ref{jfzeta} and \ref{bslemme}. This takes only a few lines. The second
section proves \ref{jfzeta}. In a third section we show that the Riemann
Hypothesis holds if and only if certain functions considered by
B\'aez-Duarte are square\--integrable, and we conclude with some related
comments.

\section{Proof of \ref{luisthmB}}

\from\  \ref{jfzeta} we have in particular on the critical line a uniform
upper estimate for $0\leq A\leq 1/2$:
$$\left|\frac{\zeta(s)}{\zeta(s+A)}\right| \leq C|s|^{1/4}$$
which implies that the functions $\zeta(s)/s\zeta(s+A)$ converge in
square-mean on the critical line to $1/s$ as $A\to0$. If we now pick one
such fixed small $A =\epsilon$ we have from \ref{bslemme} that
$$\frac{\zeta(s)}s\sum_{n=1}^N \frac{\mu(n)}{n^{s+\epsilon}} =
\frac{\zeta(s)}{s\zeta(s+\epsilon)} + O_{\epsilon,
\theta}\left(\frac{|s|^\theta}{N^{\epsilon/3}}\left|\frac{\zeta(s)}s\right|\right)$$
We then only need to invoke (as $\theta>0$ may be chosen arbitrarily)
a weak bound like $|\zeta(s)| = O(|s|^{1/4})$ (\cite[5.1.8]{titch51})
on the critical line to conclude that the left hand side converge in
square-mean sense as $N\to\infty$ to
$\zeta(s)/\zeta(s+\epsilon)s$. This concludes the proof of
\ref{luisthmB}.

It is apparent from this proof that any $\theta<1/4$ in \ref{bslemme}
will do, even any $\theta<1/2$ if we are ready to use the Lindel\"of
Hypothesis, which is a known corollary to the Riemann Hypothesis
(\cite[XIV]{titch51}). Similarly we used only a coarse version of
\ref{jfzeta}. What is essential nevertheless is the uniformity as
$A\to0$ which is granted by \ref{jfzeta}.

\section{Proof of \ref{jfzeta}}

We restate \ref{jfzeta}:

\begin{theorem}
If the Riemann Hypothesis holds one has for $\epsilon>0$
$$\left|\frac{\zeta(s)}{\zeta(s+A)}\right| =
0_\epsilon(|s|^{\inf(\epsilon,A/2)})$$
on the critical line, uniformly for $0\leq A<\infty$.
\end{theorem}

\begin{proof}
We fix $\epsilon>0$ and we may assume $\epsilon < 1/4$. The range
$2\leq A$ obviously reduces to the known Lindel\"of Hypothesis bound
$|\zeta(s)| = O(|s|^\epsilon)$.  For the range $2\epsilon\leq A \leq
2$ it is enough to combine $|\zeta(s)| = O(|s|^{\epsilon/2})$ on the
critical line with $|1/\zeta(s)| = O(|s|^{\epsilon/2})$ on $\Re(s)\geq
1/2 + 2\epsilon$ which is known to be true under the Riemann
Hypothesis (\cite[XIV]{titch51}). Finally for the range $0\leq A\leq
2\epsilon$ we invoke the next lemma.
\end{proof}

\begin{lemma}
If the Riemann Hypothesis holds one has on the critical line
$$\left|\frac{\zeta(s)}{\zeta(s+A)}\right| = O(|s|^{A/2})$$
uniformly for $0\leq A\leq1/2$.
\end{lemma}

\begin{proof}
Assuming the validity of the Riemann Hypothesis we consider the function
$$\frac{\zeta(s-A/2)}{\zeta(s+A/2)}$$
for $\Re(s)=1/2$. Let $\gamma_+(s)$ be the function appearing in the
functional equation of the zeta function:
$$\gamma_+(s) = \frac{\pi^{-s/2}\Gamma(s/2)}{\pi^{-(1-s)/2}\Gamma((1-s)/2)}$$
The following uniform estimate
$$|\gamma_+(w)| = O(|w|^{\sigma - \frac12})$$
on $\Re(w) = \sigma$,  $1/4\leq\sigma\leq3/4$, is known
(\cite[4.12.3]{titch51}).

So on the critical line and uniformly for $0\leq A\leq 1/2$:
$$\left|\frac{\zeta(s-A/2)}{\zeta(s+A/2)}\right| = |\gamma_+(s+A/2)| =
O(|s|^{A/2})$$

Next we consider for a fixed $A>0$ the holomorphic function of $s$, for
$\Re(s)\geq1/2$, given by
$$\frac{(s-1-A/2)\zeta(s-A/2)}{(s-1+A/2)\zeta(s+A/2)}\frac1{s^{A/2}}$$
On a closed vertical strip $\frac12\leq \Re(s)\leq  \Lambda$ we are in a
position to apply the Phragmen-Lindel\"of principle. Modest growth
information is necessary inside the strip for $|\zeta(s-A/2)|$ and for
$|1/\zeta(s+A/2)|$ and in fact as in the previous proof both functions are
$O(|s|^C)$ for a suitable $C$. We then let $\Lambda\to\infty$ and the
conclusion is that
$$\left|\frac{(s-1-A/2)}{(s-1+A/2)}\frac{\zeta(s-A/2)}{\zeta(s+A/2)}\right|
= O(|s|^{A/2})$$
on the closed half-plane $\Re(s)\geq1/2$, with an implied constant which is
uniform with respect to $0<A\leq 1/2$. We may of course include $A=0$ now.
The Lemma then follows from looking at the line $\Re(s) = \frac12 +
\frac{A}2$.
\end{proof}

\begin{remark}
A slightly more direct method (in place of the Phragmen-Lindel\"of
principle) is to invoke the theory of Nevanlinna functions and a result
such as \cite[Lemma 2.2]{jfjnt}. This or the simpler Phragmen-Lindel\"of
principle will also extend the bound \ref{jfzeta} from the critical line to
the closed half-plane $\Re(s)\geq1/2$ (avoiding a neighborhood of $s=1$ of
course.)
\end{remark}

\begin{remark}
The bound $|{\zeta(s-\epsilon)}/{\zeta(s+\epsilon)}| =
|\gamma_+(s+\epsilon)| = O(|s|^{\epsilon})$ on the critical line also
plays an important r\^ole in B\'aez-Duarte's proof
\cite{luisnatural}. We comment more on this proof in the next section.
\end{remark}

\section{On the functions $f_\epsilon$ considered by B\'aez-Duarte and their relation with
the Riemann Hypothesis}

B\'aez-Duarte considers for $\epsilon>0$ the function of $t>0$, defined
pointwise as:
$$f_\epsilon(t) = \sum_{n=1}^\infty \frac{\mu(n)}{n^\epsilon}\{\frac1{nt}\}$$
for $0<\epsilon$. From the bound $\{1/nt\}\leq 1/nt$ one has absolute
convergence and $|f_\epsilon(t)|\leq \zeta(1+\epsilon)/t$. Another
expression is given by:
$$f_\epsilon(t) = \frac1{\zeta(1+\epsilon) t} - \sum_{n=1}^\infty
\frac{\mu(n)}{n^\epsilon}[\frac1{nt}]$$
where the second term is pointwise a finite sum.

Let $C_n$ be the unitary operator $\phi(t)\to\sqrt{n}\phi(nt)$. One has
pointwise $f_\epsilon(t) = \sum_{n\geq1}
\mu(n)n^{-\epsilon-1/2}C_n(\{1/t\})$, so for $\epsilon>1/2$ the series is
absolutely convergent in $L^2$ and the $L^2$-Mellin transform
$\wh{f_\epsilon}(s)$ is
$$\sum_{n\geq1}
\frac{\mu(n)}{n^{\epsilon+1/2}}\frac1{n^{s-1/2}}\frac{-\zeta(s)}s =
\frac{-\zeta(s)}{\zeta(s+\epsilon)s}$$ It is known (from
\cite[3.6.5]{titch51}) that $1/\zeta(s)$ is $O(|s|^\eta)$ in
$\Re(s)\geq1$ ($\eta>0$ arbitrary) and also that $\zeta(s) =
O(|s|^{1/4})$ on $\Re(s)=1/2$ (\cite[5.1.8]{titch51}). So as
$\epsilon\to{\frac12}^+$ the functions $\wh{f_\epsilon}(s)$  on the
critical line converge in square mean. Clearly $f_\epsilon(t)\to
f_{1/2}(t)$ pointwise so this implies that $f_{1/2}(t)$ is
square-integrable. Let us recapitulate these simple results:

\begin{lemma}
For $\epsilon\geq\frac12$ the functions $f_\epsilon(t)$ are
square-integrable and their $L^2$-Mellin transforms are:
$$\wh{f_\epsilon}(s) = \frac{-1}{\zeta(s+\epsilon)}\frac{\zeta(s)}s$$
\end{lemma}

It is apparent from this expression that the square-integrability of
$f_\epsilon(t)$ for $\epsilon<1/2$ should be related with the Riemann
Hypothesis. Indeed B\'aez-Duarte proves, conditionally on the Riemann
Hypothesis, that these functions are square integrable (at least for small
$\epsilon$). He actually proves the stronger result that
$t^{-\epsilon/2}f_\epsilon(t)$ is square-integrable and that its
($L^2$)-Mellin Transform is
$-\zeta(s-\epsilon/2)/(s-\epsilon/2)\zeta(s+\epsilon/2)$. This is done so
that the limit $\epsilon\to0$ is easy to deal with. But we can copy his
argument and apply it directly to $f_\epsilon(t)$: under the Riemann
Hypothesis, and using the Theorem of Balazard-Saias \ref{bslemme}, the
($L^2$)-Mellin transforms of the finite partial sums converge in
square-mean to the function $-\zeta(s)/\zeta(s+\epsilon)s$. So the function
$f_\epsilon(t)$ which is their pointwise limit must be square-integrable.
We sum this up as a Lemma, essentially contained in \cite{luisnatural}:

\begin{lemma}
If the Riemann Hypothesis holds the functions $f_\epsilon(t)$ for
$0<\epsilon<1/2$ are square-integrable and their $L^2$-Mellin transforms are:
$$\wh{f_\epsilon}(s) = \frac{-1}{\zeta(s+\epsilon)}\frac{\zeta(s)}s$$
\end{lemma}

What we wish to point out here is the simple observation that the
square\--integrability of the function $f_\epsilon(t)$ in itself allows to
say things about the Riemann Hypothesis quite independently of its
connection with the Nyman-Beurling criterion.

Let $\epsilon>0$. We have $f_\epsilon(t) = O(1/t)$ and in particular the
function $\int_0^1 f_\epsilon(u)u^{s-1}\,du$ is analytic for $\Re(s)>1$. We
first identify this function unconditionally:

\begin{lemma}
One has unconditionally for $\Re(s)>1$:
$$\int_0^1 f_\epsilon(u)u^{s-1}\,du = \frac1{\zeta(1+\epsilon)}\frac1{s-1}
- \frac1{\zeta(s+\epsilon)}\frac{\zeta(s)}s$$
\end{lemma}

\begin{proof}
\begin{eqnarray*}
\int_0^1 f_\epsilon(u)u^{s-1}\,du
&=& \int_0^1 \left(\frac1u\frac1{\zeta(1+\epsilon)} - \sum_{a=1}^\infty
\frac{\mu(a)}{a^\epsilon}[\frac1{au}]\right)u^{s-1}\,du\\
&=& \frac1{\zeta(1+\epsilon)}\frac1{s-1} - \sum_{a=1}^\infty
\frac{\mu(a)}{a^\epsilon}\int_0^1 [\frac1{au}]u^{s-1}\,du\\
&=& \frac1{\zeta(1+\epsilon)}\frac1{s-1} - \sum_{a=1}^\infty
\frac{\mu(a)}{a^\epsilon}\int_0^{1/a} [\frac1{au}]u^{s-1}\,du\\
&=& \frac1{\zeta(1+\epsilon)}\frac1{s-1} - \sum_{a=1}^\infty
\frac{\mu(a)}{a^{s+\epsilon}}\int_0^{1} [\frac1{v}]v^{s-1}\,dv\\
&=& \frac1{\zeta(1+\epsilon)}\frac1{s-1} -
\frac1{\zeta(s+\epsilon)}\frac{\zeta(s)}s\\
\end{eqnarray*}
\end{proof}

\begin{theorem}
If $f_\epsilon\in L^2$ then its $L^2$-Mellin  transform on the
critical line is the function
$$\frac{-1}{\zeta(s+\epsilon)}\frac{\zeta(s)}s$$
Furthermore this function is analytic in $\Re(s)\geq\frac12$ except for a
simple pole at $s=1$.
\end{theorem}

\begin{proof}
Let us assume that $f_\epsilon\in L^2$. Its restriction to $t>1$ is a
non-zero multiple of $1/t$ and the corresponding Mellin transform is
$-\frac1{\zeta(1+\epsilon)}\frac1{(s-1)}$.
Let us write $k_\epsilon$ for the restriction of $f_\epsilon$ to the
interval $0<t<1$. Its Mellin transform defines an element of the Hardy
space of $\Re(s)>\frac12$, in particular it has to be analytic there. We
have obviously for $\Re(s) = \frac12$:
$$\widehat{f_\epsilon}(s) = \widehat{k_\epsilon}(s)
+\frac1{\zeta(1+\epsilon)}\frac1{1-s}$$
and we have computed explicitely $\widehat{k_\epsilon}(s)$ for $\Re(s)>1$
hence also for $\Re(s)\geq1/2$. This gives the formula
$$\frac{-1}{\zeta(s+\epsilon)}\frac{\zeta(s)}s$$ and that it
cannot have any pole  in $\Re(s)\geq1/2$ apart from $s=1$.
\end{proof}

So the square\--integrability of $f_\epsilon(t)$ is already a strong
statement; we will be content with a trivial observation:

\begin{corollary}
If $f_\epsilon\in L^2$  and $\rho$ is a zero with
$\Re(\rho)\geq\frac12+\epsilon$ then $\rho-\epsilon$ is again a zero.
\end{corollary}

This leads to:

\begin{theorem}
If $f_\epsilon\in L^2$ for a sequence of $\epsilon>0$ going to $0$ then the
Riemann Hypothesis holds.
\end{theorem}

\begin{proof}
The $\rho-\epsilon$ have to be among the zeros for an infinite sequence of
$\epsilon\to0$, but they accumulate at $\rho$.
\end{proof}

We have considered the functions $f_\epsilon(t)$ but we could as well
consider the functions $t^{-\epsilon/2}f_\epsilon(t)$ ($0<\epsilon<1/2$).
But here arises the situation that under the hypothesis that
$t^{-\epsilon/2}f_\epsilon(t)$ is square\--integrable its Fourier-Mellin
transform on the critical line has to be the function
$$\frac{-\zeta(s-\epsilon/2)}{\zeta(s+\epsilon/2)}\frac1{s-\epsilon/2}$$
which is already known to be square\--integrable!

Let us define the function  $g_\epsilon(t)$  to be the inverse
($L^2$)-Mellin transform of this function from the critical line to the
half-axis $0<t<\infty$. The Riemann Hypothesis then becomes equivalent to
the statement that $g_\epsilon(t) = t^{-\epsilon/2}f_\epsilon(t)$ (for
small $\epsilon$'s) or even seemingly weaker statements like $g_\epsilon(t)
= C_\epsilon/t^{1+\epsilon/2}$ for $t>1$. We state here just one such
strange reformulation of the Riemann Hypothesis:

\begin{theorem}
The Riemann hypothesis holds if and only if for small $\epsilon>0$ (a
sequence going to $0$ is enough), and for all $\Re(z)<\frac12$ one has:
$$\frac1{2\pi}\int_{\Re(s)=\frac12}
\frac{\zeta(s-\frac\epsilon2)}{\zeta(s+\frac\epsilon2)}\frac1{s-\frac\epsilon2}
\frac{|ds|}{s-z} = \frac1{\zeta(1+\epsilon)}\frac1{z-\frac\epsilon2 - 1}$$
\end{theorem}

\begin{proof}
Let us assume that the integral formula is true. As is known this kind of
integral corresponds to orthogonal projection to the Hardy space of the left half-plane
$\Re(z)<\frac12$. Indeed the right hand side is up to a constant the Mellin
transform of the function ${\bf 1}_{x>1}(x)x^{-1-\epsilon/2}$. But this
then implies that
$$\frac{\zeta(s-\frac\epsilon2)}{\zeta(s+\frac\epsilon2)}\frac1{s-\frac\epsilon2
}
- \frac1{\zeta(1+\epsilon)}\frac1{s-\frac\epsilon2 - 1}$$
belongs to the Hardy space of the \emph{right} half-plane $\Re(s)>\frac12$.
\from\  the absence of poles one then sees using a sequence of $\epsilon\to0$
that this implies the Riemann Hypothesis.

If the Riemann Hypothesis holds then
$$\frac{\zeta(s-\frac\epsilon2)}{\zeta(s+\frac\epsilon2)}\frac1{s-\frac\epsilon2
}$$
is the Mellin transform of $-x^{-\epsilon/2}f_\epsilon(x)$ which restricts
to $-x^{-1-\epsilon/2}/\zeta(1+\epsilon) $ for $x>1$ and this gives its
orthogonal projection to the Hardy space of $\Re(z)<\frac12$:
$$\frac{-1}{\zeta(1+\epsilon)}\int_1^\infty x^{-1-\epsilon/2} x^{z-1}\,dx =
\frac1{\zeta(1+\epsilon)}\frac1{z-\frac\epsilon2 - 1}$$
\end{proof}

We may  envision this formulation as a kind of \emph{causality}
statement about the Riemann Zeta function: the ratio
$\zeta(s-\epsilon/2)/\zeta(s+\epsilon/2)$ has the same modulus ($\Re(s)=1/2$) but not
the same phase as the spectral function $\gamma_+(s + \epsilon/2)$
of the Fourier cosine transform. In \cite{jfjnt} we formulated,
generally speaking, the Riemann Hypothesis for all abelian
$L$-functions as a statement of causality. And we expressed elsewhere
\cite{jfhab} the wish to see this implemented in a novel framework,
more truly arithmetical in nature than the ones we have been working
in in this and other papers.

\textbf{Ackowledgements.} I thank L. B\'aez-Duarte for communicating
his manuscript and also M. Balazard and \'E. Saias for the lively
discussion that followed.

\bibliographystyle{amsplain}

\end{document}